\documentclass[12pt]{article} 
\usepackage{latexsym} 
\usepackage{amsmath}
\usepackage{amssymb} 
\usepackage{amsfonts}
\usepackage{amscd} 
\usepackage{graphicx}
\usepackage{layout}
\begin{document} 
%%%%%%%%%%%%%%%%%%%%%%%%%%%%%%%%%
\newtheorem{Th}{Theorem}[section]
\newtheorem{Cor}{Corollary}[section]
\newtheorem{Prop}{Proposition}[section]
\newtheorem{Lem}{Lemma}[section]
\newtheorem{Def}{Definition}[section]
\newtheorem{Rem}{Remark}[section]
\newtheorem{Ex}{Example}[section]
\newtheorem{stw}{Proposition}[section]

%Definitions of bet ent

\newcommand{\bet}{\begin{Th}}
\newcommand{\ent}{\stepcounter{Cor}
   \stepcounter{Prop}\stepcounter{Lem}\stepcounter{Def}
   \stepcounter{Rem}\stepcounter{Ex}\end{Th}}

%Definitions of bec enc bep enp bel enl
%bef enf ber enr bee ene

\newcommand{\bec}{\begin{Cor}}
\newcommand{\enc}{\stepcounter{Th}
   \stepcounter{Prop}\stepcounter{Lem}\stepcounter{Def}
   \stepcounter{Rem}\stepcounter{Ex}\end{Cor}}
\newcommand{\bep}{\begin{Prop}}
\newcommand{\enp}{\stepcounter{Th}
   \stepcounter{Cor}\stepcounter{Lem}\stepcounter{Def}
   \stepcounter{Rem}\stepcounter{Ex}\end{Prop}}
\newcommand{\bel}{\begin{Lem}}
\newcommand{\enl}{\stepcounter{Th}
   \stepcounter{Cor}\stepcounter{Prop}\stepcounter{Def}
   \stepcounter{Rem}\stepcounter{Ex}\end{Lem}}
\newcommand{\bef}{\begin{Def}}
\newcommand{\enf}{\stepcounter{Th}
   \stepcounter{Cor}\stepcounter{Prop}\stepcounter{Lem}
   \stepcounter{Rem}\stepcounter{Ex}\end{Def}}
\newcommand{\ber}{\begin{Rem}}
\newcommand{\enr}{
   %\stepcounter{Rem} 
   \stepcounter{Th}\stepcounter{Cor}\stepcounter{Prop}
   \stepcounter{Lem}\stepcounter{Def}\stepcounter{Ex}\end{Rem}}
\newcommand{\bee}{\begin{Ex}}
\newcommand{\ene}{
 %\stepcounter{Ex}
   \stepcounter{Th}\stepcounter{Cor}\stepcounter{Prop}
   \stepcounter{Lem}\stepcounter{Def}\stepcounter{Rem}\end{Ex}}
\newcommand{\Pf}{\noindent{\it Proof\,}:\ }

%%%%%%%%%%%%%%%%%%%%%%%%%%%%%%%%%%%%%%%%%
%Beginning of Local Definition
%Local definitions
\newcommand{\EE}{\mathbf{E}}
\newcommand{\QQ}{\mathbf{Q}}
\newcommand{\R}{\mathbf{R}}
\newcommand{\C}{\mathbf{C}}
\newcommand{\ZZ}{\mathbf{Z}}
\newcommand{\KK}{\mathbf{K}}
\newcommand{\NN}{\mathbf{N}}
\newcommand{\PP}{\mathbf{P}}
\newcommand{\HH}{\mathbf{H}}
\newcommand{\uuu}{\boldsymbol{u}}
\newcommand{\xxx}{\boldsymbol{x}}
\newcommand{\aaa}{\boldsymbol{a}}
\newcommand{\bbb}{\boldsymbol{b}}
\newcommand{\AAA}{\mathbf{A}}
\newcommand{\BBB}{\mathbf{B}}
\newcommand{\ccc}{\boldsymbol{c}}
\newcommand{\ddd}{\boldsymbol{d}}
\newcommand{\iii}{\boldsymbol{i}}
\newcommand{\jjj}{\boldsymbol{j}}
\newcommand{\kkk}{\boldsymbol{k}}
\newcommand{\rrr}{\boldsymbol{r}}
\newcommand{\FFF}{\boldsymbol{F}}
\newcommand{\yyy}{\boldsymbol{y}}
\newcommand{\ppp}{\boldsymbol{p}}
\newcommand{\qqq}{\boldsymbol{q}}
\newcommand{\nnn}{\boldsymbol{n}}
\newcommand{\vvv}{\boldsymbol{v}}
\newcommand{\eee}{\boldsymbol{e}}
\newcommand{\fff}{\boldsymbol{f}}
\newcommand{\www}{\boldsymbol{w}}
\newcommand{\0}{\boldsymbol{0}}
\newcommand{\lon}{\longrightarrow}
\newcommand{\ga}{\gamma}
\newcommand{\pa}{\partial}
\newcommand{\QED}{\hfill $\Box$}
\newcommand{\id}{{\mbox {\rm id}}}
\newcommand{\Ker}{{\mbox {\rm Ker}}}
\newcommand{\grad}{{\mbox {\rm grad}}}
\newcommand{\ind}{{\mbox {\rm ind}}}
\newcommand{\rot}{{\mbox {\rm rot}}}
\newcommand{\diver}{{\mbox {\rm div}}}
\newcommand{\Gr}{{\mbox {\rm Gr}}}
\newcommand{\LC}{\mbox{\rm {LC}}}
\newcommand{\lc}{{\mbox{\rm lc}}}
\newcommand{\LG}{{\mbox {\rm LG}}}
\newcommand{\Diff}{{\mbox {\rm Diff}}}
\newcommand{\Symp}{{\mbox {\rm Symp}}}
\newcommand{\Ct}{{\mbox {\rm Ct}}}
\newcommand{\Uns}{{\mbox {\rm Uns}}}
\newcommand{\rank}{{\mbox {\rm rank}}}
\newcommand{\sign}{{\mbox {\rm sign}}}
\newcommand{\Spin}{{\mbox {\rm Spin}}}
\newcommand{\Sp}{{\mbox {\rm sp}}}
\newcommand{\Int}{{\mbox {\rm Int}}}
\newcommand{\Hom}{{\mbox {\rm Hom}}}
\newcommand{\Tan}{{\mbox {\rm Tan}}}
\newcommand{\codim}{{\mbox {\rm codim}}}
\newcommand{\ord}{{\mbox {\rm ord}}}
\newcommand{\Iso}{{\mbox {\rm Iso}}}
\newcommand{\corank}{{\mbox {\rm corank}}}
\def\mod{{\mbox {\rm mod}}}
\newcommand{\pt}{{\mbox {\rm pt}}}
\newcommand{\qed}{\hfill $\Box$ \par}
\newcommand{\spe}{\vspace{0.4truecm}}
\newcommand{\ad}{{\mbox{\rm ad}}}
\newcommand{\dint}[2]{{\displaystyle\int}_{{\hspace{-1.9truemm}}{#1}}^{#2}}

%\newcommand{\dfrac}[2]{\frac{\displaystyle{#1}}{\displaystyle{#2}}}
%%%%%%%%%%%%%%%%%%%%%%%%%%%%%%%%%%%%%%%%%
%End of Local definitions 
%%%%%%%%%%%%%%%%%%%%%%%%%%%%%%%%%

\title{
Leibniz complexity of Nash functions on differentiations
} 

\author{G. Ishikawa\thanks{This work was supported by KAKENHI No.15H03615 and 15K13431.}, 
and T. Yamashita}

%\subjclass[2010]{Primary 14P20; Secondary 58A07, 32C07}% Subject code(s)

%\keywords{K{\" a}hler differentials, complexity, extension theorem of Efroymson}% Key word(s)

%\author{G. Ishikawa\thanks{This work was supported by JSPS KAKENHI No.15H03615 and No.15K13431.}, 
%and T. Yamashita}
%
\renewcommand{\thefootnote}{\fnsymbol{footnote}}
%\footnotetext{2010 {\it Mathematics Subject Classification.} 
%Primary 14P20; Secondly 58A07, 32C07. }
%\footnotetext{{\it Key words and phrases.} K{\" a}hler differentials, complexity, extension theorem of Efroymson. 
%}
\footnotetext{The first author was supported by JSPS KAKENHI No.15H03615 and No.15K13431.}

\date{ }

\maketitle

\begin{abstract} 
The derivatives of Nash functions are Nash functions which are derived
algebraically from their minimal polynomial equations. In this paper we
show that, for any non-Nash analytic function, it is impossible to derive
its derivatives algebraically, i.e., by using linearity and Leibniz rule
finite times. In fact we prove the impossibility of such kind of algebraic computations, 
algebraically by using K{\" a}hler differentials. Then the
notion of Leibniz complexity of a Nash function is introduced in this
paper, as a computational complexity on its derivative, 
by the minimal number of usages of Leibniz rules to compute the
total differential algebraically. We provide general observations and
upper estimates on Leibniz complexity of Nash functions, related to the binary
expansions, the addition chain complexity, the non-scalar complexity and 
the complexity of Nash functions in the sense of Ramanakoraisina, 
\end{abstract}

\section{Introduction}

Let $f = f(x_1, \dots, x_n)$ be a $C^\infty$ function on an open subset $U \subset \R^n$. 
Then $f$ is called a {\it Nash function} on $U$ if $f$ is analytic-algebraic on $U$, 
i.e. if $f$ is analytic on $U$ 
and there exists a non-zero polynomial 
$P(x, y) \in \R[x, y]$, $x = (x_1, \dots, x_n)$, 
such that $P(x, f(x)) = 0$ for any $x \in U$ (\cite{Nash}\cite{Shiota}\cite{BCR}). 
If $U$ is semi-algebraic, then, 
$f$ is a Nash function if and only if $f$ is analytic and the graph of $f$ in 
$U\times\R \subset \R^{n+1}$ 
is a semi-algebraic set (\cite{BCR}). 
For a further significant progress on global study of Nash functions, see \cite{CRS04}.

An analytic function $f$ on $U$ is called {\it transcendental} if it is not a Nash function. 
Then in this paper we show that, 
for any transcendental function, it is impossible to algebraically derive 
its derivatives by using linearity and Leibniz rule (product rule) 
finite times, even by using any $C^\infty$ function. 
In fact an analytic function $f$ is a Nash function if and only if its derivatives 
$\frac{\pa f}{\pa x_1}, \dots, \frac{\pa f}{\pa x_n}$ are computable algebraically 
(Theorem \ref{characterization-theorem}). 
For example, for the transcendental function $f(x) = e^x$, the formula 
$$
\frac{d}{dx}e^x = e^x
$$ 
is never proved algebraically but is proved only by a \lq\lq transcendental" method. 
The statement above is formulated in terms of K{\" a}hler differential exactly. 

We begin with the simple example of Nash function $f(x) = \sqrt{x^2 + 1}$ of one variable. 
Then $f^2 - (x^2 + 1) = 0$. By differentiating both sides of the relation, we have 
$2f'f - 2x = 0$ where $f' = \frac{df}{dx}$. Here we have used Leibniz rule three times to get $(f^2)' = 2f'f$, $(x^2)' = 2x$ and $1' = 0$ by setting $dx/dx = 1$. 
Then we have $f'(x) = \frac{x}{f(x)} = \frac{x}{\sqrt{x^2 + 1}}$. 
If we suppose $c' = 0$ for a constant function $c$,  then the usage of Leibniz rule is counted to be twice. 

In general, let $f$ be a Nash function on $U \subset \R^n$. Then there is a non-zero polynomial 
$P(x, y) \in \R[x, y]$, $x = (x_1, \dots, x_n)$ such that $P(x, f(x)) = 0$ for any $x \in U$. 
We pose the condition that $\frac{\pa P}{\pa y}(x, f(x))$ is not identically zero on $U$. 
The condition is achieved by choosing $P$ which has the minimal total degree or the minimal degree on $y$, 
among polynomials $P$ satisfying $P(x, f(x)) = 0$ on $U$. 
Then, by using Leibniz rule in several times, we have 
$$
\frac{\pa P}{\pa x_i}(x, f(x)) \ + \ \frac{\pa P}{\pa y}(x, f(x))\frac{\pa f}{\pa x_i}(x) \ = \ 0, \quad (1 \leq i \leq n). 
$$
Therefore we have the formula 
$$
\frac{\pa f}{\pa x_i}(x) \ = \ - \ \frac{\pa P}{\pa x_i}(x, f(x))\ \big{/} \ \frac{\pa P}{\pa y}(x, f(x)), \quad 
(1 \leq i \leq n)  \quad \cdots\cdots  (*). 
$$
By our assumption that $f$ is a Nash function and the assumption on $P$, 
$\frac{\pa P}{\pa y}(x, f(x))$ is a Nash function which is not identically zero. Note that 
the above formula needs not give the value of $\frac{\pa f}{\pa x_i}(x)$ for any $x \in U$, but 
almost all $x \in U$, because $\frac{\pa P}{\pa y}(x, f(x))$ may have a zero point in $U$. 

The problem on differentiations reminds us the problem on integrations. 
Note that the partial derivatives of Nash functions are Nash functions, while 
the integrals of Nash functions need not be Nash functions. This fact was one of reason to introduce the class of elementary 
functions in classical calculus. 
For related results, say, Liouville's theorem on integrals of elementary functions, etc., refer \cite{Rosenlicht} for instance. There 
the theory of differential fields plays a significant role likewise in the present paper 
(Proofs of Lemma \ref{field-extension} and Theorem \ref{characterization-theorem}). 

Then {\it Leibniz complexity} $\LC(f)$ 
of $f$ is defined as the minimal number of usages of Leibniz rules to compute the total differential $df$ algebraically. 
The Leibniz complexity $\LC(f)$ of a Nash function $f$ is a kind of computational complexity. 
Assume any algorithm to compute the differentials of Nash functions using  
$C^\infty$ functions possibly. Then $\LC(f)$ gives 
an lower bound of usage count of Leibniz rule in such any algorithm. 
Actually we will define three variants of Leibniz complexities $\widetilde{\LC}, \LC$ and $\lc$ 
in \S \ref{Algebraic computability of differentials}. 
In particular, Nash functions are characterized by the finiteness of Leibniz complexity $\LC$ 
(Theorem \ref{characterization-theorem}). 

We remark that our complexity is closely related to the addition chain complexity \cite{Knuth} 
and to other several known computational complexities 
\cite{BS}\cite{KP}. 
We also remark that our complexity of Nash function is of different kind from the complexity for the 
description or encoding of a Nash function defined in \cite{CRS01}. 

In general it is a difficult problem to determine the exact value of the Leibniz complexity for a given Nash function. 
In \S \ref{Estimates on Leibniz complexity}, we provide general observations and 
estimates on Leibniz complexity of Nash functions using the binary expansions (Proposition \ref{estimate-2}) 
and discuss their relations with known notions on 
complexity of Nash functions (\cite{Ramanakoraisina}). 

In \S \ref{Algebraic differentiation of Nash manifolds}, we generalize Theorem \ref{characterization-theorem} to 
Nash functions on an affine Nash manifold (Theorem \ref{characterization-theorem2}), 
by using the global results on Nash functions 
(\cite{CRS96}\cite{CS00}\cite{CRS04}). 

The authors thank to anonymous referees for their valuable comments and suggestions. 
In particular the relations of Leibniz complexity with the addition chain complexity (\cite{Knuth})
and the non-scalar complexity (\cite{BS}\cite{KP}), and moreover, 
the results, Lemma \ref{monomial of one variable}, Remark \ref{lc-polynomial} 
Lemma \ref{composition-estimate} and 
Remark \ref{non-scalar complexity} are suggestions to the authors by one of the referees. 

The authors dedicate this paper to 
the memory of Professor Masahiro Shiota, who passed away  in January 2018.

\section{Algebraic computability of differentials}
\label{Algebraic computability of differentials}

Let ${\mathcal C}^{\infty}(U)$ (resp. ${\mathcal C}^{\omega}(U)$, ${\mathcal N}^{\omega}(U)$) 
denote the set of all $C^\infty$ functions (resp. analytic functions, Nash functions) on an open subset 
$U \subset \R^n$. The notation ${\mathcal N}^{\omega}(U)$ is used in \cite{Shiota}. 

Regarding $A = {\mathcal C}^{\infty}(U)$ (resp. ${\mathcal C}^{\omega}(U)$, ${\mathcal N}^{\omega}(U)$)  
as an $\R$-algebra, we take the space 
$\Omega_A$ of K{\" a}hler differentials of $A$ 
and the universal derivation ${\ddd} : A \to \Omega_A$. 

In fact, for any $\R$-algebra $A$, $\Omega_A$ can be constructed as follows: 
First consider the free $A$-module ${\mathfrak F}_A$ 
generated by elements ${\ddd}f$, for any $f \in A$, regarded as just symbols. 

Second consider the sub-$A$-module 
${\mathfrak R}_A \subset {\mathfrak F}_A$ generated by the set $R$ of 
all relations of algebraic derivations: 
$$
\ddd(h + k) - \ddd h - \ddd k, \ \ddd(\lambda \ell) - \lambda \ddd\ell, 
\ \ddd(1), \ \ddd(pq) - p\ddd q - q\ddd p, 
$$
$h, k, \ell, p, q \in A, \,  \lambda \in \R$. 
Note that an element of ${\mathfrak R}_A$ is a finite sum $\sum h_i r_i$ where $h_i \in A, \, r_i \in R$. 
Each $h_i r_i$ is called a {\it term} of the element. 
The first two kinds of generators of ${\mathfrak R}_A$ in $R$ correspond to the linearity,  
$\ddd(1)$ corresponds to the annihilation of $\R \subset A$, and 
the last kind of generators correspond to the Leibniz rule. 
We will count just the number of terms involving the last kind of generators. 
Here we add $\ddd(1)$, which 
is generated from $\ddd(1\cdot 1) - 1\ddd(1) - 1\ddd(1)$, as a generator of ${\mathfrak R}_A$ 
because we want to use the annihilation of $\R \subset A$ freely.

Third we set 
$\Omega_A = {\mathfrak F}_A/{\mathfrak R}_A$ and define 
$\ddd : 
A \to \Omega_A$ by mapping each 
$f \in A$ to the class of $\ddd f$ in ${\mathfrak F}_A/{\mathfrak R}_A$. 
Thus, 
if an element $\alpha \in {\mathfrak F}_A$ reduces to zero in $\Omega_A$, then 
there exists an element $\sum h_i r_i \in {\mathfrak R}_A$, which is called an {\it expression} of $\alpha$, 
such that $\alpha = \sum h_i r_i$ in ${\mathfrak F}_A$.

If $B$ is any $A$-module and $D : A \to B$ is any derivation, i.e. $D$ is an $\R$-linear map satisfying 
$D(g h) = gD(h) + hD(g)$ for any $g, h \in A$, then there exists a unique $A$-homomorphism 
$\rho : \Omega_A \to B$ such that $D = \rho\circ \ddd$. 

\

Suppose $U$ is connected. 

Consider the set $S \subset {\mathcal N}^{\omega}(U)$ of non-zero Nash functions i.e. 
Nash functions which are not identically zero on $U$. 
Then $S$ is closed under the multiplication. 
For $A = {\mathcal C}^{\infty}(U)$ (resp. ${\mathcal C}^{\omega}(U)$, ${\mathcal N}^{\omega}(U)$), 
let $\widetilde{A} = \widetilde{{\mathcal C}}^{\infty}(U)$ (resp. ${\widetilde{\mathcal C}}^{\omega}(U)$, 
$\widetilde{{\mathcal N}}^{\omega}(U)$) denote the localization $A_S$ of $A$ by $S$. 
Note that any element $k \in \widetilde{A}$ is 
expressed as $k = (1/g)h$ for a $g \in {\mathcal N}^{\omega}(U)$, $g \not= 0$, 
and $h \in A$ and, in general, $k$ needs not belong to $A$ if $g$ has a zero point in $U$. 
In particular ${\widetilde{\mathcal N}}^{\omega}(U) = {\mathcal N}^{\omega}(U)_S$ is 
the quotient field $Q({\mathcal N}^{\omega}(U) )$. 

Then we consider the space $\Omega_{\widetilde{A}}$ of K{\" a}hler differentials of 
the $\R$-algebra $A$ for $A = {\mathcal C}^{\infty}(U), {\mathcal C}^{\omega}(U), 
{\mathcal N}^{\omega}(U)$, 

\

Then we have: 

\bet
\label{characterization-theorem}
Let $U$ be a semi-algebraic connected open subset of $\R^n$. 
Let $A = {\mathcal C}^{\infty}(U)$ (resp. ${\mathcal C}^{\omega}(U)$, ${\mathcal N}^{\omega}(U)$). 
Then the following $10$ 
conditions on an analytic function $f \in {\mathcal C}^{\omega}(U)$ are equivalent to each other: 

{\rm (1)} $f$ is a Nash function on $U$. 

${\mbox{\rm (2)}}_A$\ 
There exists a non-zero Nash function $g \in {\mathcal N}^{\omega}(U)$ such that 
$$
g\left(\ddd f \ - \ \sum_{i=1}^n \frac{\pa f}{\pa x_i}\ddd x_i\right) = 0, 
$$
in the space $\Omega_A$ of K{\" a}hler differentials of $A$. 

${\mbox{\rm (3)}}_A$ \
$
\ddd f \ = \ \sum_{i=1}^n \frac{\pa f}{\pa x_i}\ddd x_i, 
$
in the space $\Omega_{\widetilde{A}}$ of K{\" a}hler differentials of $\widetilde{A}$. 

${\mbox{\rm (4)}}_A$ \
There exist $f_1, \dots, f_n \in \widetilde{A}$ such that 
$
\ddd f \ = \ \sum_{i=1}^n f_i\ddd x_i, 
$
in the space $\Omega_{\widetilde{A}}$ of K{\" a}hler differentials of $\widetilde{A}$. 
\ent

\

We will show the implications 
$$
\begin{array}{ccccccccc}
{\mbox{\rm (1)}} & \Rightarrow & {\mbox{\rm (2)}}_{{\mathcal N}^{\omega}(U)} & 
\Rightarrow & {\mbox{\rm (2)}}_{{\mathcal C}^{\omega}(U)} & 
\Rightarrow & {\mbox{\rm (2)}}_{{\mathcal C}^{\infty}(U)} &  & 
\vspace{0.2truecm}
\\
 & & \Downarrow & & \Downarrow & & \Downarrow & & 
\\
 &  & {\mbox{\rm (3)}}_{{\mathcal N}^{\omega}(U)} & 
\Rightarrow & {\mbox{\rm (3)}}_{{\mathcal C}^{\omega}(U)} & 
\Rightarrow & {\mbox{\rm (3)}}_{{\mathcal C}^{\infty}(U)} &  & 
\vspace{0.2truecm}
\\
& & \Downarrow & & \Downarrow & & \Downarrow & & 
\\
 &  & {\mbox{\rm (4)}}_{{\mathcal N}^{\omega}(U)} & 
\Rightarrow & {\mbox{\rm (4)}}_{{\mathcal C}^{\omega}(U)} & 
\Rightarrow & {\mbox{\rm (4)}}_{{\mathcal C}^{\infty}(U)} & \Rightarrow & {\mbox{\rm (1)}}
\end{array}
$$
to have the equivalence of the $10$ conditions. 

\

To show Theorem \ref{characterization-theorem}, we first 
recall the following known basic result on Nash functions, which is formulated in more general setting than we are going to use. 

\bel
\label{known-characterization}
Let $U \subset \R^n$ be a semi-algebraic open subset and 
$f \in {\mathcal C}^{\omega}(U)$ be an analytic function on $U$. 
Then the following conditions are equivalent to each other:

{\rm (i)}
$f$ is a Nash function on $U$, i.e. 
there exists a non-zero polynomial $P(x, y)$ such that $P(x, f(x)) = 0$ for any $x \in U$. 

{\rm (ii)}
The graph of $f$ in $U\times\R \subset \R^{n+1}$ 
is a semi-algebraic set. 

{\rm (iii)}
For any $a \in U$, the Taylor series $j^\infty f(a)$ of $f$ at $a$ is algebraic in the 
of formal power series algebra $\R[[x - a]]$  over 
the polynomial algebra $\R[x - a] = \R[x]$, in other words, there exists a non-zero polynomial 
$P(x, y)$ such that $j^\infty P(x, f)(a) = 0$. 

{\rm (iv)}
For any connected component $U'$ of $U$,  
there exists a point $a \in U'$ such that the Taylor series $j^\infty f(a)$ of $f$ at $a$ is algebraic in 
formal power series algebra $\R[[x - a]]$  over 
the polynomial algebra $\R[x - a]$. 
\enl

\Pf
The equivalences (i) and (ii) are well-known (see for instant \cite{BCR}). 
The implications (i) $\Rightarrow$ (iii) 
$\Rightarrow$ (iv) are clear. To show the implication (iv) $\Rightarrow$ (i), 
suppose (iv). Note that the number of connected components of $U$ is finite. 
Let $U_1, \dots, U_r$ are all connected components of $U$. Let $1 \leq i \leq r$. 
Then there exists $a_i \in U_i$ such that 
$f$ is expressed by the Taylor series at $a_i$ in a neighborhood $W \subset U_i$ of $a_i$ and 
there exists a non-zero polynomial $P_i(x, y)$ such that $P_i(x, f(x)) = 0$ for any $x \in W$. 
Since the function $P_i(x, f(x))$ is analytic on $U_i$ and $U_i$ is connected, $P_i(x, f(x)) = 0$ for any $x \in U_i$. 
Then it suffices to take $P = \prod_{i=1}^r P_i$ to get (i). 
\QED

\

Also we need the general algebraic lemma to show the implication 
${\mbox {\rm (4)}}_{{{\mathcal C}^\infty}(U)} \Rightarrow$ (1) of 
Theorem \ref{characterization-theorem}. 

\bel
\label{field-extension}
Let $K \subset L$ be a field extension. Assume that $\R \subset K$. 
Let $f \in L$ be a transcendental element over $K$. 
Then, for any derivation $D_0 : K \to K$ and for any $u \in L$, 
there exists a unique derivation $D_u : K(f) \to L$ satisfying 
$$
D_u\vert_K = D_0, \quad D_u(f) = u. 
$$
Moreover if $L$ is finitely generated over $K$, 
then the derivation $D_u$ extends to a derivation $D : L \to L$. 
\enl

\Pf
Since $f$ is transcendental over $K$, 
we can define a derivation $D_u : K(f) \to L$ on the extension field $K(f)$ over $K$ by $f$, 
by $D_u\vert_K = D_0$ and $D_u(f) = u$. 
Suppose $L$ is finitely generated over $K$ and $L = K(f, h_1, \dots, h_m)$ for some 
$h_1, \dots, h_m \in L$. 
Then we define a derivation $D_{u1} : K(f, h_1) \to L$, $D_{u1}\vert_{K(f)} = D_u$ as follows: 
If $h_1$ is transcendental over $K(f)$, then we set $D_{u1}(h_1) = 0$. If 
$h_1$ is algebraic over $K(f)$, then we set $D_{u1}(h_1)$ as the element in $K(f, h_1)$ 
which is determined by the algebraic relation of $h_1$ over $K(f)$ and $D_u$. 
In fact, if 
$\sum_{k=0}^m a_k h_1^{m-k} = 0, \ a_k \in K(f), $ is 
a minimal algebraic relation of $h_1$ over $K(f)$, then we would have 
$$
\sum_{k=0}^m D_u(a_k) h_1^{m-k} + \left(\sum_{k=0}^{m-1} (m-k) a_k h_1^{m-k-1}\right) D_1(h_1)= 0. 
$$
Since $\sum_{k=0}^{m-1} (m-k) a_k h_1^{m-k-1} \not= 0$ by the minimality assumption, 
$D_{u1}(h_1)$ is uniquely determined by 
$$
D_{u1}(h_1) = - \left(\sum_{k=0}^m D_u(a_k) h_1^{m-k}\right) \big{/} \left(\sum_{k=0}^{m-1} (m-k) a_k h_1^{m-k-1}\right). 
$$
Thus we extend $D_u$ into a derivation $D = D_{um} : L \to L$ by a finitely number of steps. 
Note that we need not to use Zorn's lemma to show the existence of extension of derivation. 
\QED

\

\noindent
{\it Proof of Theorem \ref{characterization-theorem}.} 
(1) $\Rightarrow$ ${\mbox{\rm (2)}}_{{\mathcal N}^{\omega}(U)}$ : 
Let $f \in {\mathcal C}^{\omega}(U)$ be a Nash function and $P(x, y)$ be a non-zero polynomial satisfying 
$P(x, f) = 0$ and $\frac{\pa P}{\pa y}(x, f) \not= 0$. 
Then, by taking K{\" a}hler differential on both sides of the polynomial equality $P(x, f) = 0$, 
we have in $\Omega_{{\mathcal N}^{\omega}(U)}$, 
\begin{eqnarray*}
0 & = & \ddd(P(x, f)) = \sum_{i=1}^n \frac{\pa P}{\pa x_i}(x, f)\ddd x_i + 
\frac{\pa P}{\pa y}(x, f)\ddd f 
\\
& = & \sum_{i=1}^n \left(- \frac{\pa P}{\pa y}(x, f)\frac{\pa f}{\pa x_i}\right)\ddd x_i 
+ \frac{\pa P}{\pa y}(x, f)\ddd f 
\\
& = & \frac{\pa P}{\pa y}(x, f)\left( \ddd f - \sum_{i=1}^n \frac{\pa f}{\pa x_i}\ddd x_i \right), 
\end{eqnarray*}
and that $\frac{\pa P}{\pa y}(x, f)$ is a non-zero Nash function on $U$. 

Since ${\mathcal N}^{\omega}(U) \subset {\mathcal C}^{\omega}(U) \subset 
{\mathcal C}^{\infty}(U)$, 
the implications 
${\mbox{\rm (j)}}_{{\mathcal N}^{\omega}(U)}$ $\Rightarrow$ 
${\mbox{\rm (j)}}_{{\mathcal C}^{\omega}(U)}$ $\Rightarrow$ 
${\mbox{\rm (j)}}_{{\mathcal C}^{\infty}(U)}$ are clear, for ${\mbox{\rm j}} = 2, 3, 4$. 

${\mbox{\rm (2)}}_A$ $\Rightarrow$ ${\mbox{\rm (3)}}_A$, $A =  {\mathcal C}^{\infty}(U), 
{\mathcal C}^{\omega}(U), {\mathcal N}^{\omega}(U)$: 
Since $1/g$ belongs to the localization ${\widetilde{A}}$, we have that, 
if $g\left(\ddd f \ - \ \sum_{i=1}^n \frac{\pa f}{\pa x_i}\ddd x_i\right) = 0$ in 
$\Omega_A$, then 
$\ddd f \ - \ \sum_{i=1}^n \frac{\pa f}{\pa x_i}\ddd x_i = 0$ in 
$\Omega_{\widetilde{A}}$. 

%The implications 
%${\mbox{\rm (3)}}_{{\mathcal N}^{\omega}(U)}$ $\Rightarrow$ 
%${\mbox{\rm (3)}}_{{\mathcal C}^{\omega}(U)}$ $\Rightarrow$ 
%${\mbox{\rm (3)}}_{{\mathcal C}^{\infty}(U)}$ are clear. 

The implications ${\mbox{\rm (3)}}_A$ $\Rightarrow$ ${\mbox{\rm (4)}}_A$, $A =  {\mathcal C}^{\infty}(U), 
{\mathcal C}^{\omega}(U), {\mathcal N}^{\omega}(U)$, are clear. 

%The implications 
%${\mbox{\rm (4)}}_{{\mathcal N}^{\omega}(U)}$ $\Rightarrow$ 
%${\mbox{\rm (4)}}_{{\mathcal C}^{\omega}(U)}$ $\Rightarrow$ 
%${\mbox{\rm (4)}}_{{\mathcal C}^{\infty}(U)}$ are clear. 

${\mbox{\rm (4)}}_{{\mathcal C}^{\infty}(U)}$ $\Rightarrow$ (1) : 
Suppose $f$ is not a Nash function on $U$ and 
$\ddd f - \sum_{i=1}^n f_i\ddd x_i = 0$ in $\Omega_{\widetilde{\mathcal C}^{\infty}(U)}$. 
Since $f$ is not a Nash function, by Lemma \ref{known-characterization}, 
there exists a point $a \in U$ such that $f \in \R[[x - a]] 
\subset 
Q(\R[[x - a]])$ is not algebraic. Here 
$\R[[x - a]]  = {\mathcal C}^{\infty}_{\R^n, a}/{\mathfrak m}^\infty_{\R^n, a}$ is the $\R$-algebra of 
formal series, $M = Q(\R[[x - a]])$ is its quotient field and 
the Taylor series of $f$ at $a$ is written also by the same symbol $f$. 
Moreover, we have $\ddd f - \sum_{i=1}^n f_i\ddd x_i = 0$ 
in the K{\" a}hler differentials $\Omega_M$ of $M$, 
via the homomorphism $\widetilde{\mathcal C}^{\infty}(U) \to M$ defined by taking the Taylor series. 
Then, in the free $M$-module ${\mathfrak F}_M$ 
generated by elements 
$\{\ddd h \mid h \in M\}$, 
$\ddd f - \sum_{i=1}^n f_i\ddd x_i$ is a finite sum of elements of type 
$$
a(\ddd(h + k) - \ddd h - \ddd k), \ b(\ddd(\lambda \ell) - \lambda \ddd\ell), \ c(\ddd(pq) - p\ddd q - q\ddd p).
$$
Here $a, h, k, b, \ell, c, p, q \in M, \lambda \in \R$. 
Now we take the subfield $L \subset M$ generated over the rational function field $K = \R(x)$ by 
$f, f_i (1 \leq i \leq n)$ 
and those $a, h, k, b, \ell, c, p, q$ which appear in the above expression of 
$\ddd f - \sum_{i=1}^n f_i\ddd x_i$: 
$L = K(f, h_1, \dots, h_m)$, which is a finitely generated field over $K$ by $f$ and for some 
$h_1, \dots, h_m \in M$. 
Then we have 
$\ddd f - \sum_{i=1}^n f_i\ddd x_i = 0$ also in $\Omega_{L}$. 

Take any non-zero element $u \in L$ and fix it. Set $D_0 = 0$. 
Then, by Lemma \ref{field-extension}, we have a derivation $D : L \to L$ with $D(f) = u$. 
Then by the universality of 
the K{\" a}hler differentials, there exists an $L$-linear map
$\rho : \Omega_{L} \to L$ such that 
$\rho\circ \ddd = D : L \to L$. Here $\ddd : L \to \Omega_L$ is the universal derivation. 
Then we have 
$$
0 = \rho\left(\ddd f - \sum_{i=1}^n f_i\ddd x_i\right) = D(f) = u. 
$$ 
This leads to a contradiction with the assumption $u \not= 0$. 
Thus we have that $f$ is a Nash function.  
\qed

\ber
{\rm 
If Zorn's lemma is used, 
then the fact that a transcendental basis of $M = Q(\R[[x - a]])$ forms a basis of $\Omega_M$ as an  
$M$-vector space (Theorem 26.5 of \cite{Matsumura}) will give a shorter proof of 
the part ${\mbox{\rm (4)}}_{{\mathcal C}^{\infty}(U)}$ $\Rightarrow$ (1) 
of proof of Theorem \ref{characterization-theorem}. In fact 
if $f \in M$ is transcendental, then there exists a transcendental basis containing $f, x_1, \dots, x_n$ and 
therefore we have that $\ddd f, \ddd x_1, \dots, \ddd x_n$ are linearly independent over $M$, which leads a contradiction. 
(The remark is based on an anonymous reviewer's comment informed to the authors.) \ The same remark 
is applied also to the proof of our Theorem \ref{characterization-theorem2}. 
}
\enr

\ber
{\rm
If $U$ is not connected, then Theorem \ref{characterization-theorem}
does not hold.
In fact, let $U = \R \setminus \{ 0\}$ and set $f(x) = e^x$ if $x > 0$ and
$f(x) = 1$ if $x < 0$.
Then $f \in {\mathcal C}^{\omega}(U)$ and $f \not\in {\mathcal N}^{\omega}(U)$.
However the condition (2) is satisfied if we take as $g$ the non-zero Nash
function
on $U$ defined by $g(x) = 0 (x > 0), \, g(x) = 1 (x < 0)$.
}
\enr

\section{Estimates on Leibniz complexity}
\label{Estimates on Leibniz complexity}

Let $U \subset \R^n$ be a semi-algebraic connected open subset. 
Let $f \in {\mathcal N}^{\omega}(U)$ be a Nash function on $U$. 

Then by the equivalence of (1) and ${\mbox{\rm (2)}}_{{\mathcal N}^{\omega}(U)}$ in 
Theorem \ref{characterization-theorem}, 
there exists a non-zero Nash function $g \in {\mathcal N}^{\omega}(U)$ such that 
$g(\ddd f - \sum_{i=1}^n \frac{\pa f}{\pa x_i}\ddd x_i) 
\in {\mathfrak R}_{{\mathcal N}^{\omega}(U)} (\subset {\mathfrak F}_{{\mathcal N}^{\omega}(U)}$). 
Then define $\LC_g(f)$ as the minimal number of 
terms corresponding to Leibniz rule for all expressions of 
$g(\ddd f - \sum_{i=1}^n \frac{\pa f}{\pa x_i}\ddd x_i) 
\in {\mathfrak R}_{{\mathcal N}^{\omega}(U)}$. 
We define the {\it Leibniz complexity} ${\mbox{\rm LC}}(f)$ of $f$ by the minimum of 
$\LC_g(f)$ for all such non-zero $g \in {\mathcal N}^{\omega}(U)$. 

Note that we do not care about the number of terms corresponding to linearity of the differential. 
Moreover we do not count the term generated by the relation $\ddd(1\cdot 1) - 1\ddd(1) - 1\ddd(1)$. 
Therefore we use the relation $\ddd(c) = 0$ for $c \in \R$ freely. 

Similarly we define $\widetilde{\LC}(f)$, related to Theorem 
\ref{characterization-theorem} ${\mbox{\rm (2)}}_{{\mathcal N}^{\omega}(U)}$, 
as the minimal number of 
terms corresponding to Leibniz rule for all expressions of 
$\ddd f - \sum_{i=1}^n \frac{\pa f}{\pa x_i}\ddd x_i
\in {\mathfrak R}_{{\widetilde{\mathcal N}}^{\omega}(U)}$. 

Moreover, if $\ddd f - \sum_{i=1}^n \frac{\pa f}{\pa x_i}\ddd x_i \in {\mathfrak R}_{{\mathcal N}^{\omega}(U)}$, 
we define $\lc(f) = \LC_1(f)$,  
simply as the minimal number of 
terms corresponding to Leibniz rule for all its expressions in ${\mathfrak R}_{{\mathcal N}^{\omega}(U)}$. 
Note that, if $f$ is a polynomial function, then $\lc(f) < \infty$. 
However in general $\ddd f - \sum_{i=1}^n \frac{\pa f}{\pa x_i}\ddd x_i$ may not 
belong to ${\mathfrak R}_{{\mathcal N}^{\omega}(U)}$. Then we set $\lc(f) = \infty$. 

\

Hereafter, for $f \in {\mathcal N}^{\omega}(U)$, 
we set 
$$
\zeta(f) := \ddd f - \sum_{i=1}^n \frac{\pa f}{\pa x_i}\ddd x_i, 
$$
regarded as an element in ${\mathfrak F}_A$ for $A = {\mathcal N}^\omega(U)$ or 
for its localization 
$A = {\widetilde{\mathcal N}}^\omega(U) = {\mathcal N}^\omega(U)_S$ 
where $S = {\mathcal N}^\omega(U) \setminus \{ 0\}$. 

First we show general basic inequalities:

\bel
\label{inequalities for three invariants}
For any $f \in {\mathcal N}^{\omega}(U)$, we have \ 
$
\widetilde{\LC}(f) \ \leq \ \LC(f) \ \leq \ \lc(f). 
$
\enl

\Pf
Suppose $\lc(f) < \infty$ and there exists an expression of $\zeta(f)$ 
in ${\mathfrak R}_{{\mathcal N}^\omega(U)}$ 
such that the number of terms involving Leibniz rule is equal to $\lc(f)$. 
Then setting $g = 1$, $g\zeta(f)$ has the same expression in ${\mathfrak R}_{{\mathcal N}^\omega(U)}$, 
and therefore we have $\LC(f) \leq \lc(f)$. 
Next, by the definition of  ${\LC}(f)$, there exist 
a $g \in {\mathcal N}^\omega(U) \setminus \{ 0\}$ and 
an expression of $g\zeta(f)$ in ${\mathfrak R}_{{\mathcal N}^\omega(U)}$ such that 
the number of terms involving Leibniz rule is equal to ${\LC}(f)$. 
Then, dividing by $g$, we have an expression of $\zeta(f)$ in 
${\mathfrak R}_{{\widetilde{\mathcal N}}^\omega(U)}$ such that 
the number of terms involving Leibniz rule is equal to ${\LC}(f)$. 
Therefore, by the definition of $\widetilde{\LC}(f)$, we have $\widetilde{\LC}(f) \leq {\LC}(f)$. 
\QED

\bel
For $f, g \in {\mathcal N}^{\omega}(U)$, we have 

{\rm (1)} 
$
\LC(f + g) \leq \LC(f) + \LC(g). 
$
\quad 
{\rm (2)} 
$
\LC(fg) \leq \LC(f) + \LC(g) + 1. 
$

The same inequalities hold for $\widetilde{\LC}$ and $\lc$. 
\enl

\Pf
Let $h\zeta(f) \in {\mathfrak R}_{{\mathcal N}^{\omega}(U)}$ (resp. $k\zeta(g) 
\in {\mathfrak R}_{{\mathcal N}^{\omega}(U)}$) 
be expressed using the terms of Leibniz rule minimally i.e. $\LC(f)$-times (resp. 
$\LC(g)$-times), for a non-zero $h \in {\mathcal N}^{\omega}(U)$ (resp. 
a non-zero $k \in {\mathcal N}^{\omega}(U)$).
Then $hk\zeta(f + g) = k(h\zeta (f) ) +h(k\zeta (g) ) \in {\mathfrak R}_{{\mathcal N}^{\omega}(U)}$ 
is expressed using 
Leibniz rule at most $\LC(f) + \LC(g)$ times. Therefore we have (1). Moreover, 
by using Leibniz rule once, we have 
$$
hk\ddd(fg) = hk(g\ddd f + f\ddd g) = kg(h\ddd f) + hf(k\ddd g)
$$ 
in $\Omega_{{\mathcal N}^{\omega}(U)}$.  
Then, using Leibniz rule $\LC(f) + \LC(g)$ times, we compute $h\ddd f$ 
and $k\ddd g$, and thus $hk\ddd(fg)$. Therefore we have (2). 

For $\widetilde{\LC}$ and $\lc$, the inequalities are proved similarly or more easily. 
\QED

\

By the definition of Leibniz complexity, we have the affine invariance: 

\bel
\label{affine-isom}
Let $f \in {\mathcal N}^{\omega}(U)$ and $\varphi : \R^n \to \R^n$ be an affine isomorphism. 
Then $f\circ\varphi \in {\mathcal N}^{\omega}(\varphi^{-1}(U))$ 
satisfies $\LC(f\circ\varphi) = \LC(f)$, $\widetilde{\LC}(f\circ\varphi) = \widetilde{\LC}(f)$ 
and $\lc(f\circ\varphi) = \lc(f)$. 
\enl

\Pf
By the definition of Leibniz complexity 
$h(\ddd f - \sum_{i=1}^n \frac{\pa f}{\pa x_i}\ddd(x_i))$ is zero in $\Omega_{{\mathcal N}^{\omega}(U)}$ 
by using Leibniz rule $\LC(f)$-times, for a non-zero $h \in {\mathcal N}^{\omega}(U)$. 
Let $x' = (x_1', \dots, x_n')$ be new affine coordinate system on $\R^n$ defined by $x' = \varphi^{-1}(x)$. 
Then
$(h\circ\varphi)(\ddd (f\circ\varphi) - \sum_{i=1}^n \frac{\pa f}{\pa x_i}\circ\varphi \ddd(\varphi_i))$ is zero 
in $\Omega_{{\mathcal N}^{\omega}(U)}$ by using Leibniz rule $\LC(f)$-times.
Since we do not count the usage of Leibniz rule for $\ddd(c)= 0, c \in \R$, we have that 
$(h\circ\varphi)(\ddd (f\circ\varphi) - \sum_{i=1}^n \frac{\pa f\circ\varphi}{\pa x'_i} \ddd(x'_i))$ is 
zero in $\Omega_{{\mathcal N}^{\omega}(U)}$ by using Leibniz rule the same $\LC(f)$-times. 
Note that $h\circ\varphi \in {\mathcal N}^{\omega}(U)$ is non-zero. Therefore we have $\LC(f\circ\varphi) \leq \LC(f)$. 
Similarly, we have $\LC(f) = \LC((f\circ\varphi)\circ\varphi^{-1}) \leq \LC(f\circ\varphi)$. Thus we have the required equality. 
The equality for $\widetilde{\LC}$ (resp. $\lc(f)$) is proved similarly or more easily. 
\QED

\

In general it is a difficult problem to determine the exact value of the Leibniz complexity even for an  
polynomial function. 

\bee
{\rm
Let $n = 1$ and write $x = x_1$. Then we have 
$\widetilde{\LC}(x + c) = \LC(x + c) = \lc(x + c) = 0$. 
$\widetilde{\LC}(x^2 + bx + c) = \LC(x^2 + bx + c) =  \lc(x^2 + bx + c) = 1$. 
$\widetilde{\LC}(\sqrt{x^2 + 1}) = \LC(\sqrt{x^2 + 1}) = \lc(\sqrt{x^2 + 1}) = 2$. 

Let $n = 2$.  For $\lambda \in \R$, we have
$$
{\mathrm {LC}}(x_1^2 + x_2^2 + \lambda x_1x_2) = 
\left\{ 
\begin{array}{ccc}
1 \ & {\mbox{\rm if}} & \vert \lambda\vert \geq 2 
\\
2 \ & {\mbox{\rm if}} & \vert \lambda\vert < 2. 
\end{array}
\right.
$$
In fact, 
$x_1^2 + x_2^2 + \lambda x_1x_2 = (x_1 + \frac{\lambda}{2}x_2)^2 + (1 - \frac{\lambda^2}{4})x_2^2$. 
Moreover $x_1^2 + x_2^2 + \lambda x_1x_2 = (x_1 + \alpha x_2)(x_1 + \beta x_2)$ 
for some $\alpha, \beta \in \R$
if and only if $\vert \lambda\vert \geq 2$. The same results hold for $\widetilde{\LC}$ and $\lc$. 
}
\ene

\

Let $n = 1$ and write $x = x_1$. We consider Leibniz complexity of a monomial $x^k$. 
For example, $\lc(x^0) = \lc(1) = 0, \lc(x) = 0, \lc(x^2) = 1, \lc(x^3) = 2, \lc(x^4) = 2$. 
Also for $\LC$ and $\widetilde{\LC}$ we have the same results. 
For example we calculate $\ddd(x^4) = 2x^2\ddd(x^2) = 4x^3\ddd (x)$ by using Leibniz rule 
twice, and we can check that it is impossible to calculate $\ddd(x^4)$ by using Leibniz rule 
just once.  

\

To observe the essence of the problem to estimate the Leibniz complexity, 
let us digress to consider \lq\lq the problem of strips\rq\rq. Let 
$k$ be a positive integer. Suppose we have a sheet of paper having width $k$ 
and, using a pair of scissors, we make $k$-strips of width $1$. 
We may cut several sheets of the same width 
at once by piling them. Then the problem is to minimize the total number of cuts. 
Clearly it is at most $k-1$. 

The exact answer to the above problem is given by the {\it addition chain complexity} $\ell(k)$ (see \cite{Knuth}). 
An {\it addition chain} of $k$ is a sequence of integers 
$$
1 = a_0, \ a_1, \ a_2, \ \dots, \ a_r = k 
$$
satisfying that, for any $i = 1, 2, \dots, r$, there exists $j, m$ with $0 \leq j \leq m < i$, 
such that $a_i = a_j + a_m$. 
Then $\ell(k)$ is defined as the minimum of the length $r$ for all addition chain of $k$. 

A process of making $k$-strips as above corresponds to 
an addition chain bijectively. Therefore the minimum of the total number of cuts is 
given by $\ell(k)$. 

\bel
\label{monomial of one variable}
For a positive integer $k$, we have 
$$
\widetilde{\LC}(x^k) \leq \LC(x^k) \leq \lc(x^k) \leq \ell(k). 
$$
\enl

\Pf
Let $1 = a_0, \ a_1, \ a_2, \ \dots, \ a_r = k$ be an addition chain of $k$. Since 
$k = a_r = a_j + a_m$ for some $0 \leq j \leq m < k$, we have one relation
$$
\ddd(x^k) - x^{a_j}\ddd(x^{a_m}) - x^{a_m}\ddd(x^{a_j})
$$
in ${\mathfrak R}_{{\mathcal N}^\omega(U)}$. Thus we have
$$
\ddd(x^k) = x^{a_j}\ddd(x^{a_m}) + x^{a_m}\ddd(x^{a_j})
$$
in $\Omega_{{\mathcal N}^\omega(U)}$ using Leibniz rule once. 
If $j < m$, we apply this procedure to $\ddd(x^{a_m})$. Then, 
using a relation 
$$
\ddd(x^k) - x^{a_j}\ddd(x^{a_m}) - x^{a_m}\ddd(x^{a_j}) 
+ x^{a_j}(\ddd(x^{a_m}) - x^{a_{j'}}\ddd(x^{a_{m'}}) - x^{a_{m'}}\ddd(x^{a_{j'}}))
$$
with two terms, in the sense of \S \ref{Algebraic computability of differentials}, 
in ${\mathfrak R}_{{\mathcal N}^\omega(U)}$ for some $0 \leq j' \leq m' < m$, 
we have 
$$
\ddd(x^k) = x^{a_j + a_{j'}}\ddd(x^{a_{m'}}) + x^{a_j + a_{m'}}\ddd(x^{a_{j'}})
$$
in $\Omega_{{\mathcal N}^\omega(U)}$ using Leibniz rule twice. 
If $j = m$, then $x^{a_j}\ddd(x^{a_m}) + x^{a_m}\ddd(x^{a_j}) = 2x^{a_m}\ddd(x^{a_m})$, 
and then a similar procedure is applied to $\ddd(x^{a_m})$. 
Thus we see that, by using a relation with $s$-terms involving Leibniz rule, 
$\ddd(x^k)$ is reduced to a functional linear combination of 
$\ddd(x^{a_0}), \ddd(x^{a_1}), \dots, \ddd(x^{a_{r - s}})$ in $\Omega_{{\mathcal N}^\omega(U)}$, 
$s = 1, 2, \dots, r$. 
Therefore we have $\lc(x^k) (= \LC_1(x^k)) \leq r$, for any addition chain of $k$. 
Hence we have $\lc(x^k) \leq \ell(k)$. 
Other inequalities follow from Lemma \ref{inequalities for three invariants}. 
\QED

\ber
\label{lc-polynomial}
{\rm 
We can define, naturally, a kind of Leibniz complexity $\lc_{\mbox{\rm {\scriptsize poly}}}$ 
by using the K{\" a}hler differential $\Omega_A$ of 
polynomial algebra 
$A = \R[x_1, \dots, x_n]$. Then the proof of Lemma \ref{monomial of one variable} gives also 
the inequalities $\lc(x^k) \leq \lc_{\mbox{\rm {\scriptsize poly}}}(x^k) \leq \ell(k)$. 
The authors conjecture, at least, the equality $\lc_{\mbox{\rm {\scriptsize poly}}} (x^k) = \ell(k)$, 
but they have no proof of that. 
}
\enr

Now we show one known strategy to obtain an explicit estimate. 
Consider the binary expansion of $k$: 
$$
k = 2^{\mu_r} + 2^{\mu_{r-1}} + \cdots + 2^{\mu_1},
$$
for some integers $\mu_r > \mu_{r-1} > \cdots > \mu_1 \geq 0$. We set $\mu = \mu_r$. Then 
{\it the number of digits} (\lq$1$\rq\ or \lq$0$\rq)  
is given by $\mu + 1$, while $r$ is {\it the number of units}, \lq$1$\rq, appearing in the binary expansion. 
Then first we cut the sheet into $r$ sheets of width $2^{\mu}, 2^{\mu_{r-1}}, \dots, 2^{\mu_1}$ 
by $(r-1)$-cuts. Second, divide the sheet of width $2^{\mu}$ into sheets of width $2^{\mu_{r-1}}$ 
by $\mu - \mu_{r-1}$-cuts. Third, divide the piled sheets of width $2^{\mu_{r-1}}$ into sheets of width $2^{\mu_{r-2}}$ 
by $\mu_{r-1} - \mu_{r-2}$-cuts, and so on. Iterating the process, we have sheets of width $2^{\mu_1}$, 
which we divide into strips of width $1$ by $\mu_1$-cuts finally. 
The total number of cuts by this method is given by $\mu + r - 1$. 

Thus we have by Lemma \ref{monomial of one variable}: 

\bec
\label{estimate-of-monomial}
For a positive integer $k$, we have 
$$
\widetilde{\LC}(x^k) \leq \LC(x^k) \leq \lc(x^k) \leq \ell(k) \leq \mu + r - 1. 
$$
\enc

\ber
{\rm 
The estimate in Corollary \ref{estimate-of-monomial} is, by no means, best possible. 
For example, let $k = 31$. Then $31 = 2^4 + 2^3 + 2^2 + 2^1 + 2^0$. 
Therefore $r = 5$ and $\mu = 4$. Therefore $\mu + r - 1 = 8$. 
Moreover we have the addition chain complexity $\ell(31) = 7$. 
%So the estimate gives us that $\LC(x^{31}) \leq 8$. 
However
$\LC(x^{31}) \leq 6$. In fact, since $32 =2^5$, we have by Lemma \ref{estimate-of-monomial}, 
$$
x\ddd(x^{31}) = \ddd(x^{32}) - x^{31}\ddd(x) = 32x^{31}\ddd(x) - x^{31}\ddd(x) = 31x^{31}\ddd(x), 
$$
by using Leibniz rule $6$ times. Then we have $\ddd(x^{31}) = 31x^{30}\ddd(x)$ 
in $\Omega_{{\widetilde{\mathcal C}}^{\infty}(U)}$. 
}
\enr

Related to Corollary \ref{estimate-of-monomial}, we observe 

\bel
For $f \in {\mathcal N}^{\omega}(U)$ and a natural number $k \geq 1$, we have 
$
\LC(f^k) \leq \LC(f) + \LC(x^k). 
$
\enl

\Pf
If $f$ is a constant function, then $\LC(f^k) = 0$, so the inequality holds trivially. 
We suppose $f$ is not a constant function. 
By definition, for some non-zero $g \in {\mathcal N}^{\omega}(\R)$, 
$g\ddd(x^k)$ is deformed into $g\,kx^{k-1}\ddd x$ in 
$\Omega_{{\mathcal N}^{\omega}(\R)}$ using Leibniz rules $\LC(x^k)$-times. 
Using the same procedure, $(g\circ f)\ddd(f^k)$ is deformed into $(g\circ f)kf^{k-1}\ddd f$ in 
$\Omega_{{\mathcal N}^{\omega}(U)}$ using Leibniz rules $\LC(x^k)$-times. 
Note that $g\circ f$ is non-zero in ${\mathcal N}^{\omega}(U)$. 
Moreover, using Leibniz rules $\LC(f)$ times, $h(g\circ f)kf^{k-1}\ddd f$ is deformed into 
$h(g\circ f)\sum_{i=1}^n kf^{k-1}(\pa f/\pa x_i)\ddd x_i$ for some non-zero $h \in {\mathcal N}^{\omega}(U)$. 
Since $g\circ f$ is non-zero, $h(g\circ f)$ is non-zero. 
\QED

\

In general we have 

\bel
\label{composition-estimate}
Let $g_1, \dots, g_m \in {\mathcal N}^\omega(U)$ and 
$P(y_1, \dots, y_m) \in \R[y_1, \dots, y_m]$ be a polynomial 
regarded as a function on $\R^m$. 
Then, for the Leibniz complexity of $f = P(g_1, \dots, g_m)$, we have 
$$
\LC(f) \leq \lc(P) + \sum_{i=1}^m \LC(g_i), \ 
\widetilde{\LC}(f) \leq \lc(P) + \sum_{i=1}^m \widetilde{\LC}(g_i), 
\ 
\lc(f) \leq \lc(P) + \sum_{i=1}^m \lc(g_i). 
$$
\enl

\Pf
We give a proof of the first inequality only. The remaining inequalities are proved similarly or more easily. 

Using Leibniz rule $\lc(P)$ times, we have 
$$
\ddd(f) \ = \ \sum \frac{\pa P}{\pa y_i}(g_1, \dots, g_m)\,\ddd(g_i), 
$$
in $\Omega_{{\mathcal N}^\omega(U)}$. 
For each $i = 1, \dots, m$, 
there exists non-zero Nash function $h_i$ such that 
$$
h_i\ddd(g_i) = h_i\sum_{j=1}^n \frac{\pa g_i}{\pa x_j}\ddd(x_j)
$$
by an $\LC(g_i)$ times usage of Leibniz rule. Therefore 
$$
h_1\cdots h_m\ddd(f) = h_1\cdots h_m(\sum_{j=1}^n \frac{\pa f}{\pa x_j}\ddd(x_j)), 
$$
in $\Omega_{{\mathcal N}^\omega(U)}$, using Leibniz rule $\lc(P) + \sum_{i=1}^m \LC(g_i)$ times 
in total. Therefore we have $\LC(f) \leq \lc(P) + \sum_{i=1}^m \LC(g_i)$. 
\QED

\ber
\label{non-scalar complexity}
{\rm 
The Leibniz complexity $\lc(P)$ or $\lc_{\mbox{\rm {\scriptsize poly}}}(P)$ 
(see Remark \ref{lc-polynomial}) for polynomials $P$ is closely related to the {\it non-scalar complexity} 
of $P$ (\cite{KP}\cite{BS}). The non-scalar complexity of a polynomial $P$ is defined roughly as follows. 
Consider any program to produce polynomials in $\R[x_1, \dots, x_n]$ by scalar multiplications, 
additions and products, without divisions, starting 
from the $0$-th stage $1, x_1, \dots, x_n$ (depth $0$), and making some pair of linear combinations 
of polynomials appeared in previous stages of depth $\leq r$
and, as the next stage, making the product of them (depth $r + 1$) and so on. 
Then the non-scalar complexity ${\mbox{\rm L}}_{\mbox{\rm \scriptsize{ns}}}(P)$ is defined as 
the minimal depth of the polynomial $P$ in all such 
programs producing $P$. Then we have 
$$
\lc(P) \ \leq \ \lc_{\mbox{\rm {\scriptsize poly}}}(P) \ \leq \ {\mbox{\rm L}}_{\mbox{\rm \scriptsize{ns}}}(P). 
$$
The proof is similar to that of Lemma \ref{monomial of one variable}. 
The authors conjecture also that 
the equality $\lc_{\mbox{\rm {\scriptsize poly}}}(P) = {\mbox{\rm L}}_{\mbox{\rm \scriptsize{ns}}}(P)$ holds, 
but they have no proof of the equality. 

In \cite{BS}, the non-scalar complexity of rational functions for programs allowing divisions is considered 
and, for any rational function $f$, 
an estimate of the non-scalar complexity of partial derivatives $\frac{\pa f}{\pa x_i}$ 
by means of that of $f$. It is interesting to estimate the Leibniz complexity of 
partial derivatives of higher order by Baur-Strassen's result \cite{BS}. 
}
\enr

\

As above, we consider \lq\lq the problem of strips\rq\rq\ starting from several number of sheets, say, $s$,  
having width $k_s$, $k_{s-1}$, and $k_1$ respectively. Then we have 

\

\bel
\label{estimate-1}
Let $P = P(x) = a_{s}x^{k_s} + a_{s-1}x^{k_{s-1}} +\cdots + a_1x^{k_1} \in \R[x]$ 
be a polynomial function of one variable, where 
$a_j \not= 0\, (1 \leq j \leq s)$ and $k_s > k_{s-1} > \cdots > k_1 \geq 0$. 
Regarding the binary expansion, let $\mu$ be (the number of digits of $k_s$) $- 1$,  and 
$r_j$ 
the number of units of $k_j$, $1 \leq j\leq s$. 
Then, by using Leibniz rule $\mu + \sum_{j = 1}^s (r_j - 1)$-times and linearity, and by supposing $\ddd(c) = 0, c \in \R$, we have 
$\ddd(P) = (dP(x)/dx)\ddd(x)$ in $\Omega_{{\mathcal N}^{\omega}(U)}$. 
In particular we have 
$$
\widetilde{\LC}(P) \leq \LC(P) \leq \lc(P) \leq \mu + \sum_{j = 1}^s (r_j - 1). 
$$
\enl

\Pf
Let $\mu = \mu_t > \mu_{t-1} > \cdots > \mu_1 \geq 0$ be 
all of the exponents appearing in the binary expansions of $k_s, k_{s-1}, \dots, k_1$. 
First, by using Leibniz rule $\sum_{j = 1}^s (r_j - 1)$-times, we modify $\ddd(P)$ 
into a linear combination of $\ddd(x^\ell), \ell = 2^\mu = 2^{\mu_t}, 2^{\mu_{t-1}}, \dots, 2^{\mu_1}$. 
Second, by using Leibniz rule $\mu - \mu_{t-1}$-times, we modify $\ddd(x^\ell), \ell = 2^\mu$ into 
$\ddd(x^{\ell'}), \ell' = 2^{\mu_{t-1}}$. Repeating the procedure, we modify $\ddd(P)$ into 
a multiple of $\ddd(x^\ell), \ell = 2^{\mu_1}$. Finally, by using Leibniz rule $\mu_1$-times, 
we modify $\ddd(P)$ into a multiple of $\ddd(x)$. 
\QED

We estimate the Leibniz complexity for a polynomial of $n$-variables. 
Let $P(x) = P(x_1, \dots, x_n) \in \R[x_1, \dots, x_n]$. 
We set $P(x) = \sum b_\alpha x^\alpha, b_\alpha \in \R,$ by 
using multi-index $\alpha = (\alpha_1, \dots, \alpha_n)$ of non-negative integers. 
It is trivial that $\lc(P)$ is at most the total number of multiplications of variables: 
$$
\sum_{b_\alpha \not= 0} \max\{ \vert \alpha\vert - 1, \ 0\}. 
$$
Instead we consider the number 
$$
\sigma(P) := \sum_{b_\alpha \not= 0} \max\{ \#\{ i \mid 1 \leq i \leq n, \alpha_i > 0\} - 1, \ 0\}, 
$$
which is needed just to separate the variables on differentiation, and we try to save the additional usage 
of Leibniz rule. 

Suppose that, by arranging terms with respect to $x_i$ 
for each $i, 1 \leq i \leq n$, 
$$
P(x) = a_{i,s(i)}x_i^{k_{i,s(i)}} + a_{i,s(i)-1}x_i^{k_{i,s(i)-1}} +\cdots + a_{i,1}x_i^{k_{i,1}}, 
$$ 
where $a_{i,j}$ is a non-zero polynomial of $x_1, \dots, x_n$ without $x_i$, 
$(1 \leq j \leq s(i))$, and $k_{i,s(i)} > k_{i,s(i)-1} > \cdots > k_{i,1} \geq 0$. 
The maximal exponent $k_{i,s(i)}$ is written as $\deg_{x_i}P$,
the degree of $P$ in the variable $x_i$.
For the binary expansion of $\deg_{x_i}P$,
let $\mu_i$ denote (the number of digits of $\deg_{x_i}P$) $- 1$.
Moreover let $r_{ij}, 1 \leq j \leq s(i)$ denote the number of units of
the exponent $k_{ij}$
for the binary expansion.
Then we have 

\bel
\label{estimate-2}
By using the linearly, $\ddd(c) = 0, c \in \R$, and Leibniz rule 
\\
$\sigma(P) + \sum_{i=1}^n \left(\mu_i + \sum_{j = 1}^{s(i)} (r_{ij} - 1)\right)$-times, 
we have 
$\ddd(P) = \sum_{i=1}^n (\pa P(x)/\pa x_i)\ddd(x_i)$ in $\Omega_{{\mathcal N}^{\omega}(U)}$. 
In particular we have the estimate 
$$
\lc(P) \leq \sigma(P) + \sum_{i=1}^n \left(\mu_i + \sum_{j = 1}^{s(i)} (r_{ij} - 1)\right). 
$$
\enl

\ber
{\rm 
We have, for any polynomial $P(x) = \sum b_\alpha x^\alpha$, 
$$
\sigma(P) + \sum_{i=1}^n \left(\mu_i + \sum_{j = 1}^{s(i)} (r_{ij} - 1)\right)
\ \leq \ \sum_{b_\alpha \not= 0} \max\{ \vert \alpha\vert - 1, 0\}. 
$$
and in almost cases the inequality is strict. 
}
\enr

\noindent
{\it Proof of Lemma \ref{estimate-2}.} 

By applying Leibniz rule to each term of $P$, 
$\ddd(P)$ is deformed into a sum of forms $a_{i,j}\ddd(x_i^{k_{i,j}})$ 
with the differential of one variable $x_i$ and a function $a_{i,j}$ of other variables. 
For this process we need to use Leibniz rule $\sigma(P)$-times. 
Then $\ddd(P)$ is the sum of the form 
$$
a_{i,s(i)}\ddd(x_i^{k_{i,s(i)}}) + a_{i,s(i)-1}\ddd(x_i^{k_{i,s(i)-1}}) +\cdots + a_{i,1}\ddd(x_i^{k_{i,1}}), 
$$
($i = 1, \dots, n$). By Lemma \ref{estimate-1}, for each $i = 1, \dots, n$, 
the form is deformed into 
$\sum_{i=1}^n \frac{\pa P}{\pa x_i}\ddd x_i$ by using Leibniz rule $\mu_i + \sum_{j = 1}^{s(i)} (r_{ij} - 1)$. 
Thus we have the estimate. 
\QED

\

Now we give an upper estimate of Leibniz complexities for Nash functions by 
those for polynomial functions in terms of its polynomial relation. 
Let $f \in {\mathcal N}^{\omega}(U)$ be a Nash function on a connected open subset $U$ of $\R^n$. 
Let $P(x, y) = P(x_1, \dots, x_n, y)$ be a polynomial such that 
$P(x, f(x)) = 0$ on $U$ and $\frac{\pa P}{\pa y}(x, f(x))$ is not identically zero. 
We set $x_0 = y$. Suppose that, by arranging with respect to $x_i$ 
for each $i, 0 \leq i \leq n$, 
$$
P(x, y) = a_{i,s(i)}x_i^{k_{i,s(i)}} + a_{i,s(i)-1}x_i^{k_{i,s(i)-1}} +\cdots + a_{i,1}x_i^{k_{i,1}}, 
$$ 
where $a_{i,j}$ is a non-zero polynomial of $x_0, x_1, \dots, x_n$ without $x_i$, 
$(1 \leq j \leq s(i))$, and $k_{i,s(i)} > k_{i,s(i)-1} > \cdots > k_{i,1} \geq 0$. 
For the binary expansion, 
let $\mu_i$ (resp. $r_{ij}, 1 \leq j \leq s(i)$)
be (the number of digits of $\deg_{x_i}P$) $-1$
(resp. the number of units of $k_{ij}$), $0 \leq i \leq n$, respectively. 
Write $\deg_{x_i}P$ the degree of $P$ with respect to $x_i, 0 \leq i \leq n$ 
and use the same notation $\sigma(P)$ as in Lemma \ref{estimate-2} for the polynomial $P$ of 
$n+1$ variables. 

\bep
\label{estimate-3}
Under the above notations, we have the estimate 
$$
\LC(f) \leq \sigma(P) + \sum_{i=0}^n \left(\mu_i + \sum_{j = 1}^{s(i)} (r_{ij} - 1)\right). 
$$
In particular we have 
$$
\LC(f) \leq \sigma(P) + \sum_{i=0}^n\{ (\deg_{x_i}P +2)(\log_2(\deg_{x_i}P) - 1) \} + n + 1.
$$
\enp

\bee
{\rm 
Let $n = 1, f = \frac{1}{\sqrt{x^2 + 1}}$ and $P(x, y) = y^2 - x^2 - 1$. Then $\sigma(P) = 0, 
\mu_0 = \mu_1 = 1$ and $r_{ij} = 1$. Therefore the first inequality gives us 
that $\LC(f) \leq 2$ as is seen in Introduction. 
}
\ene

\noindent
{\it Proof of Proposition \ref{estimate-3}.} 

We write the right hand side by $\psi$ of the first inequality. 
By Lemma \ref{estimate-2}, we have, by using Leibniz rule $\psi$-times, 
$$
\ddd(P(x, y)) = \sum_{i=1}^n \frac{\pa P}{\pa x_i}(x, y)\ddd x_i + \frac{\pa P}{\pa y}(x, y)\ddd y, 
$$
modulo several linearity relations and $\ddd c, c \in \R$ in ${\Omega}_{{\mathcal N}^{\omega}(U\times \R)}$. 
Then, substituting $y$ by $f$, we have that 
$$
0 = \ddd(P(x, f)) = \sum_{i=1}^n \frac{\pa P}{\pa x_i}(x, f)\ddd x_i + \frac{\pa P}{\pa y}(x, f)\ddd f, 
$$
in $\Omega_{{\mathcal N}^{\omega}(U)}$, therefore that
$$
\frac{\pa P}{\pa y}(x, f)\left( \ddd f - \sum_{i=1}^n \frac{\pa f}{\pa x_i} \ddd x_i\right) = 0, 
$$
in ${\Omega}_{{\mathcal N}^{\omega}(U)}$, 
by using Leibniz rule at most $\psi$-times. Thus we have the first inequality. 
The second equality is obtained from the first equality combined with the inequalities derived by the definitions: 
$$
2^{\mu_i} \leq \deg_{x_i} P  < 2^{\mu_i+1}, \ s(i) \leq \deg_{x_i} P  + 1, {\mbox{\rm \ and }} \ r_{ij} \leq \mu_i, 
$$
$(1 \leq j \leq s(i), 0 \leq i \leq n)$. 
\QED

\

In \cite{Ramanakoraisina}, the complexity ${\mathrm{C}}(f)$ of a Nash function $f$ is defined as the minimum 
the total degree $\deg P$ of non-zero polynomials $P(x, y)$ with $P(x, f) = 0$. 
Moreover we define 
$$
{\mathrm{S}}(f) := \min\{ \sigma(P\circ \psi) \mid P(x, f) = 0, \deg P = {\mathrm{C}}(f), 
\psi {\mbox{\rm \ is an affine isomorphism on\ }} \R^{n+1} \}, 
$$
i.e. the minimum of the number $\sigma$ for any defining polynomial $P$ of $f$ with minimal total degree 
under any choice of affine coordinates. 
We can regard ${\mathrm{S}}(f)$ a complexity for the separation of variables in differentiation of $f$. 
Then we have the following result: 

\bec
Let $f \in {\mathcal N}^{\omega}(U)$ be a Nash function on a connected open set $U \subset \R^n$. 
Then we have an estimate on the Leibniz complexity $\LC(f)$ by 
the Ramanakoraisina's complexity ${\mathrm{C}}(f)$ and another complexity ${\mathrm{S}}(f)$, 
$$
\LC(f) \leq {\mathrm{S}}(f) + (n+1)({\mathrm C}(f) + 2)(\log_2\!{\mathrm C}(f) - 1) + n+1. 
$$
\enc

\Pf
Since 
$\deg_{x_i}P \leq {\mathrm C}(f) \, (0 \leq i \leq n)$ we have the above estimate by Proposition \ref{estimate-3}
and Lemma \ref{affine-isom}. 
\QED

\

Naturally we would like to pose a problem to obtain any lower estimate of Leibniz complexity. 

\section{Algebraic differentiation on Nash manifolds} 
\label{Algebraic differentiation of Nash manifolds}

Let $U$ be a connected semi-algebraic open subset of $\R^n$ and $M \subset U$ 
a Nash submanifold (\cite{BCR}\cite{Shiota}). Suppose $M$ is a closed connected subset in $U$. 
We consider the quotient $\R$-algebra ${\mathcal N}^{\omega}(U)/I$ by the ideal $I$ of ${\mathcal N}^{\omega}(U)$ consisting of Nash functions on $U$ which vanish on $M$. 

Since ${\mathcal N}^{\omega}(U)$ is Noetherian (\cite{Risler}\cite{Mostowski}), 
$I$ is generated by a finite number of Nash functions 
$g_1, \dots, g_\ell \in {\mathcal N}^{\omega}(U)$ over ${\mathcal N}^{\omega}(U)$. 

An element $[f] \in {\mathcal C}^{\omega}(U)/I{\mathcal C}^{\omega}(U)$ is called {\it Nash} if 
there exists a polynomial $P(x, y) = a_m(x)y^m + a_{m-1}(x)y^{m-1} + \cdots + a_1(x)y + a_0(x) \in \R[x, y]$ 
satisfying that at least one of $a_m([x]), a_{m-1}([x]), \dots, a_1([x]), a_0([x])$ is not zero in ${\mathcal N}^{\omega}(U)/I$ 
and that $P([x], [f]) = 0$ in ${\mathcal C}^{\omega}(U)/I{\mathcal C}^{\omega}(U)$. 
The condition is equivalent to that $[f]$ is algebraic over $\R(x)$ via the composition 
$\R(x) \hookrightarrow {\mathcal C}^{\omega}(U) \to {\mathcal C}^{\omega}(U)/I{\mathcal C}^{\omega}(U)$ of natural homomorphisms. 
Also the condition is equivalent to that $[f]$ is algebraic over ${\mathcal N}^{\omega}(U)/I$ via the natural homomorphism
${\mathcal N}^{\omega}(U)/I \to {\mathcal C}^{\omega}(U)/I{\mathcal C}^{\omega}(U)$. 
Then there exist a non-zero polynomial $P(x, y)$ and 
$h_j \in {\mathcal C}^{\omega}(U), 1 \leq j \leq \ell$ such that 
$$
P(x, f(x)) = \sum_{j=1}^\ell h_j(x)g_j(x), 
$$ 
for any $x \in U$ and that $\frac{\pa P}{\pa y}(x, f) \not\in I{\mathcal C}^{\omega}(U)$. 
By differentiating both sides of the relation by $x_i$, we have that 
$$
\frac{\pa P}{\pa x_i}(x, f(x)) + \frac{\pa P}{\pa y}(x, f(x))\frac{\pa f}{\pa x_i}
= \sum_{j=1}^\ell g_j(x) \frac{\pa h_j}{\pa x_i}(x) +  \sum_{j=1}^\ell h_j(x) \frac{\pa g_j}{\pa x_i}(x), 
$$
so that 
$$
\frac{\pa P}{\pa y}([x], [f])\left[\frac{\pa f}{\pa x_i}\right] = - \frac{\pa P}{\pa x_i}([x], [f]), 
$$ 
in ${\mathcal C}^{\infty}(U)/(I + \langle \pa g_1/\pa x_i, \dots, \pa g_\ell/\pa x_i\rangle_{{\mathcal C}^{\infty}(U)})$, 
for $1 \leq i \leq n$. Note that $\frac{\pa P}{\pa y}([x], [f])$ is non-null in ${\mathcal C}^{\omega}(U)/I{\mathcal C}^{\omega}(U)$ 
and algebraic over ${\mathcal N}^{\omega}(U)/I$. 

\

We consider the space 
$\Omega_A$ of K{\" a}hler differentials of  
%$A = {\mathcal C}^{\infty}(U)/I{\mathcal C}^{\infty}(U)$ (resp.  
%${\mathcal C}^{\omega}(U)/I{\mathcal C}^{\omega}(U), {\mathcal N}^{\omega}(U)/I$). 
$A = {\mathcal C}^{\infty}(U)$ (resp.  
${\mathcal C}^{\omega}(U), {\mathcal N}^{\omega}(U)$). 
Note that 
$
\Omega_{A/IA} \cong \Omega_{A}/(AdI + I\Omega_{A}), 
$
as an $A/IA$-module. 
For the set $S$ of non-zero Nash elements in ${\mathcal C}^{\omega}(U)/I{\mathcal C}^{\omega}(U)$, 
$\widetilde{A / I A} = (A / I A)_S$ denote the localization of 
$A / I A = {\mathcal C}^{\infty}(U)/I{\mathcal C}^{\infty}(U)$ (resp.  
${\mathcal C}^{\omega}(U)/I{\mathcal C}^{\omega}(U), {\mathcal N}^{\omega}(U)/I$) by $S$.

An ideal $I$ of ${\mathcal N}^{\omega}(U)$ is called {\it locally formally prime} if, 
for each $a \in U$, 
the ideal $I_a$ in the formal algebra $\R[[x - a]]$ generated by $\{j^\infty h(a) \mid h \in I\}$ 
is prime.

Then we have: 

\bet
\label{characterization-theorem2}
Let $U$ be a connected semi-algebraic open subset of $\R^n$ and 
$I$ a locally formally prime ideal in ${\mathcal N}^{\omega}(U)$. 
Let $A = {\mathcal C}^{\infty}(U)/I{\mathcal C}^{\infty}(U)$, 
${\mathcal C}^{\omega}(U)/I{\mathcal C}^{\omega}(U)$ or ${\mathcal N}^{\omega}(U)/I$. 
Then the following $10$ conditions on 
$[f] \in {\mathcal C}^{\omega}(U)/I{\mathcal C}^{\omega}(U)$ are equivalent to each other: 

{\rm (1)} 
$[f]$ is Nash. 

${\mbox{\rm (2)}}_A$ \ 
There exists a non-zero Nash element $[g] \in {\mathcal C}^{\omega}(U)/I{\mathcal C}^{\omega}(U)$ such that 
$$
[g]\left(\ddd [f] \ - \ \sum_{i=1}^n \left[ \frac{\pa f}{\pa x_i}\right] \ddd [x_i]\right) = 0, 
$$
in the space $\Omega_A$ 
of K{\" a}hler differentials of $A$. 

${\mbox{\rm (3)}}_A$ \ 
$
\ddd [f] \ = \ \sum_{i=1}^n \left[ \frac{\pa f}{\pa x_i}\right] \ddd [x_i], 
$
in the space $\Omega_{\widetilde{A}}$ of K{\" a}hler differentials of 
the localization $\widetilde{A}$ of $A$ by the set of non-zero Nash elements. 

${\mbox{\rm (4)}}_A$ \
There exist $\alpha_1, \dots, \alpha_n \in \widetilde{A}$ such that 
$
\ddd [f] \ = \ \sum_{i=1}^n \alpha_i \ddd [x_i], 
$
in the space $\Omega_{\widetilde{A}}$.  
\ent

\ber
{\rm 
If $I$ is the ideal of Nash functions 
vanishing on a connected closed Nash submanifold $M \subset U$, then 
$I$ is locally formally prime and 
$I{\mathcal C}^\omega(U)$ is prime in ${\mathcal C}^\omega(U)$. 
}
\enr

To show Theorem \ref{characterization-theorem2}, we need 
the following characterization of Nash function. 
It is proved using the extension theorem due to Efroymson or its 
generalization \cite{CS00}:

\bel
\label{known-characterization2}
Let $U \subset \R^n$ be a connected semi-algebraic open subset and 
$I \subset {\mathcal N}^{\omega}(U)$ be an ideal. 
For any $f \in {\mathcal C}^{\omega}(U)$ the following conditions are equivalent to each other:

{\rm (i)}
$[f] \in {\mathcal C}^{\omega}(U)/I{\mathcal C}^{\omega}(U)$ is Nash. 

{\rm (ii)}
For any $a \in U$, the Taylor series $j^\infty f(a)$ of $f$ at $a$ is algebraic in 
$\R[[x - a]]/I_a$, 
in other words, there exists a polynomial $P(x, y) \in \R[x, y], \deg_yP > 0$, which possibly depends on $a$, 
such that $j^\infty P(x, f)(a) \in I_a$, where 
$I_a$ is the ideal in $\R[[x - a]]$ generated by $\{j^\infty h(a) \mid h \in I\}$. 

{\rm (iii)} 
There exists a Nash function $g \in {\mathcal N}^{\omega}(U)$ such that $[g] = [f] \in {\mathcal C}^{\omega}(U)/I{\mathcal C}^{\omega}(U)$. 
\enl

\Pf
The implication (i) $\Rightarrow$ (ii) is clear. 

(ii) $\Rightarrow$ (iii): 
Let ${\mathcal I}$ be the finite ideal sheaf generated by $I$ in the sheaf ${\mathcal N}^{\omega}_U$ of Nash functions. 
Then $f$ defines a section of the quotient sheaf ${\mathcal N}^{\omega}_U/{\mathcal I}$. By the extension theorem (\cite{CRS96}\cite{CS00}) 
in non-compact case, there exists $g \in {\mathcal N}^{\omega}(U)$ which defines the same section of ${\mathcal N}^{\omega}_U/{\mathcal I}$ with that defined by $f$. 
Therefore $f - g  \in {\mathcal C}^{\omega}(U)$ defines a section of ${\mathcal I}{\mathcal C}^{\omega}_U$, the ideal sheaf generated by ${\mathcal I}$ in the sheaf ${\mathcal C}^{\omega}_U$ of analytic functions.  Then $f - g \in I{\mathcal C}^{\omega}(U)$, by Cartan's theorem A for real analytic functions (\cite{Cartan}). 
Thus we have (iii). 

The implication (iii) $\Rightarrow$ (i) is clear. 
\QED

\

\noindent
{\it Proof of Theorem \ref{known-characterization2}.} 
(1) $\Rightarrow$ ${\mbox{\rm (2)}}_{{\mathcal N}^{\omega}(U)/I}$: 
Suppose (1). We take a representative $f$ which belongs to ${\mathcal N}^{\omega}(U)$ 
by Lemma \ref{known-characterization2}. 
Then we have 
\begin{eqnarray*}
0 & = & \ddd(P([x], [f])) = \sum_{i=1}^n \frac{\pa P}{\pa x_i}([x], [f])\ddd [x_i] + 
\frac{\pa P}{\pa y}([x], [f])\ddd [f] 
\\
& = & \sum_{i=1}^n \left(- \frac{\pa P}{\pa y}([x], [f])\left[\frac{\pa f}{\pa x_i}\right]\right)\ddd [x_i] 
+ \frac{\pa P}{\pa y}([x], [f])\ddd [f] 
\\
& = &
\frac{\pa P}{\pa y}([x], [f])\left(\ddd [f] - \sum_{i=1}^n \frac{\pa P}{\pa x_i}([x], [f])\ddd [x_i] \right), 
\end{eqnarray*}
in $\Omega_{{\mathcal N}^{\omega}(U)/I}$, and $\frac{\pa P}{\pa y}([x], [f]) 
\in {\mathcal N}^{\omega}(U)/I$ 
is non-zero and algebraic over ${\mathcal N}^{\omega}(U)/I$. 

The implications ${\mbox{\rm (j)}}_{{\mathcal N}^{\omega}(U)/I}$ 
$\Rightarrow$ ${\mbox{\rm (j)}}_{{\mathcal C}^{\omega}(U)/I{\mathcal C}^{\omega}(U)}$ 
$\Rightarrow$ ${\mbox{\rm (j)}}_{{\mathcal C}^{\infty}(U)/I{\mathcal C}^{\infty}(U)}$ 
are clear, for ${\mbox{\rm j}} = 2, 3, 4$. 

The implications ${\mbox{\rm (2)}}_A$ $\Rightarrow$ ${\mbox{\rm (3)}}_A$, 
for $A = {\mathcal N}^{\omega}(U)/I, {\mathcal C}^{\omega}(U)/I{\mathcal C}^{\omega}(U), 
{\mathcal C}^{\infty}(U)/I{\mathcal C}^{\infty}(U)$, 
 are clear, since $[g] \in S$. 
 
%The implications ${\mbox{\rm (3)}}_{{\mathcal N}^{\omega}(U)/I}$ 
%$\Rightarrow$ ${\mbox{\rm (3)}}_{{\mathcal C}^{\omega}(U)/I{\mathcal C}^{\omega}(U)}$ 
%$\Rightarrow$ ${\mbox{\rm (3)}}_{{\mathcal C}^{\infty}(U)/I{\mathcal C}^{\infty}(U)}$ 
%are clear. 
 
The implications ${\mbox{\rm (3)}}_A$  $\Rightarrow$ ${\mbox{\rm (4)}}_A$ for 
$A = {\mathcal N}^{\omega}(U)/I, {\mathcal C}^{\omega}(U)/I{\mathcal C}^{\omega}(U), 
{\mathcal C}^{\infty}(U)/I{\mathcal C}^{\infty}(U)$ are clear. 

%The implications ${\mbox{\rm (4)}}_{{\mathcal N}^{\omega}(U)/I}$ 
%$\Rightarrow$ ${\mbox{\rm (4)}}_{{\mathcal C}^{\omega}(U)/I{\mathcal C}^{\omega}(U)}$ 
%$\Rightarrow$ ${\mbox{\rm (4)}}_{{\mathcal C}^{\infty}(U)/I{\mathcal C}^{\infty}(U)}$ 
%are clear. 

${\mbox{\rm (4)}}_{{\mathcal C}^{\infty}(U)/I{\mathcal C}^{\infty}(U)}$ $\Rightarrow$ (1): 
Suppose ${\mbox{\rm (4)}}_{{\mathcal C}^{\infty}(U)/I{\mathcal C}^{\infty}(U)}$ and $[f]$ is not Nash. 
Then, by Lemma \ref{known-characterization2}, there exists a point $a \in U$ such that 
$[f]$ is transcendental in $\R[[x - a]]/I_a$ 
via the $\R$-algebra homomorphism $\varphi_a : {\mathcal N}^{\omega}(U)/I \to \R[[x - a]]/I_a$, where 
$I_a$ is the ideal in the formal power series ring $\R[[x - a]]$ generated by $g_1, \dots, g_\ell$. 
Let $K = Q(\varphi_a({\mathcal N}^{\omega}(U)/I))$ be the quotient field of the image of ${\mathcal N}^{\omega}(U)/I$ by $\varphi_a$. 
Moreover let $L = K([f], [h_1], \dots, [h_m])$ be the extended field of $K$ 
which is generated by all elements which appear in the relation 
$\ddd [f] - \sum_{i=1}^n \alpha_i \ddd [x_i] = 0$ in $\Omega_{\R[[x - a]]/I_a}$. 
Then the relation holds also in $\Omega_L$. 

Let $u$ be any non-zero element of $L$. 
We extend the zero derivation $D_0 = 0 : K \to L$ to $D_u : K([f]) \to L$ by setting $D_u([f]) = u$, 
for the given non-zero element $u \in L$. 
Moreover we extend $D_u$ to a derivation $D : L \to L$. 
Then for an $L$-homomorphism $\rho : \Omega_L \to L$ we have $D = \rho \circ \ddd : L \to L$. 
Then we have 
$$
0 = \rho\left(\ddd [f] - \sum_{i=1}^n \alpha_i \ddd [x_i]\right) = D([f]) = u. 
$$ 
This leads a contradiction. Thus we have (1). 
\QED

\

For a Nash element $[f] \in {\mathcal C}^{\omega}(U) / I {\mathcal C}^{\omega}(U)$, 
we define the Leibniz complexity of $[f]$ by the minimal number of 
terms corresponding to Leibniz rule for $[g] \left( \ddd [f] - \sum_{i=1}^n \left[ \frac{\pa f}{\pa x_i} \right] \ddd [x_i] \right)$ 
in the free ${\mathcal N}^{\omega}(U) /I$-module 
${\mathfrak F}_{{\mathcal N}^{\omega}(U)/I}$
among all expressions for all non-zero Nash element 
$[g] \in {\mathcal C}^{\omega}(U) /  I {\mathcal C}^{\omega}(U)$. 
The definition is based on the statement 
${\mbox{\rm (2)}}_{{\mathcal N}^{\omega}(U)/I}$ of Theorem \ref{characterization-theorem2}. 
We do not care about the number of terms corresponding to linearity of the differential. 
Moreover we will do not count the term generated by the relation $\ddd([1\cdot 1]) - [1]\ddd([1]) - [1]\ddd([1])$. 
Therefore we use the relation $\ddd([c]) = 0$ for $c \in \R$ freely. 

Let ${\mathrm {LC}}([f])$ denote the Leibniz complexity of $[f]$. 
Similarly to Proposition \ref{estimate-3} we have an upper estimate: 

\bep
\label{estimate-4}
Under the situation of Theorem \ref{characterization-theorem2}, 
let $P(x, y)$ be a polynomial such that $P(x, f) \in I{\mathcal C}^{\omega}(U)$ 
and $\frac{\pa P}{\pa y}(x, f) \not\in I{\mathcal C}^{\omega}(U)$. 
Then we have 
$$
\LC([f]) \leq \sigma(P) + \sum_{i=0}^n \left(\mu_i + \sum_{j = 1}^{s(i)} (r_{ij} - 1)\right). 
$$
\enp

\

\

\begin{flushleft}
Goo ISHIKAWA, \\
Department of Mathematics, Hokkaido University, 
Sapporo 060-0810, Japan. \\
e-mail : ishikawa@math.sci.hokudai.ac.jp \\

\

Tatsuya YAMASHITA, \\
Department of Mathematics, Hokkaido University, 
Sapporo 060-0810, Japan. \\
e-mail : tatsuya-y@math.sci.hokudai.ac.jp \\ 

\end{flushleft}

\end{document}